\documentclass[reqno]{amsart}
\newtheorem{theorem}{Theorem}[]

\theoremstyle{definition}
\newtheorem{definition}[]{Definition}
\newtheorem{example}[]{Example}

\theoremstyle{remark}

\begin{document}

\title{Solution of Partial Differential Equations by Method of Hyperholomorphic functions}

\author{Anatoliy A. Pogorui}
\address{Department of Mathematics, Zhytomyr State University, Zhytomyr, Ukraine 10008}
\email{pogor@zu.edu.ua}



\date{November 16, 2008}



\begin{abstract}
It is well known that the real and imaginary parts of any
holomorphic function are harmonic functions of two variables. In
this paper we generalize this property to finite-dimensional
commutative algebras. We prove that if some basis of a subspace of
a commutative algebra satisfies a polynomial equation then
components of a hyperholomorphic function on the subspace are
solutions of the respective partial differential equation.
\end{abstract}

\maketitle

Let  $\mathbf{A}$  be a finite-dimensional commutative unitary
algebra over $K=\mathbb{R}$ (or $\mathbb{C}$), a set of vectors
$\vec{e}_{0}, \vec{e}_{1}, \ldots, \vec{e}_{n}$ be a basis of
$\mathbf{A}$, and $\vec{e}_{0}$ be the unit of the algebra.
Suppose $\mathbf{B}$ is $m$-dimensional subspace of $\mathbf{A}$
with the basis $\vec{e}_{0}, \vec{e}_{1}, \ldots, \vec{e}_{m}$,
$m\leq n$.
 Consider a function
$\vec{f}:\mathbf{B}\rightarrow\mathbf{A}$ of the following form
$$
 \vec{f}(\vec{x})=\sum_{k=0}^{n}\vec{e}_{k}u_{k}(\vec{x}),
$$
where $u_{k}(\vec{x})=u_{k}(x_{0},x_{1},\ldots,x_{m})$  are real
(or complex) functions of $m+1$ arguments.

\begin{definition}
$\vec{f}(\vec{x})$ is called differentiable at a point
$\vec{x}_{0}\in \mathbf{B}$ if there exists the function
$\vec{f'}:\mathbf{B}\rightarrow\mathbf{A}$  such that for any
$\vec{h}\in \mathbf{B}$
\begin{equation}\label{1}
\vec{h}\vec{f'}(\vec{x}_{0})=\lim_{\varepsilon\rightarrow 0}
\frac{\vec{f}(\vec{x}_{0}+\varepsilon\vec{h})-\vec{f}(\vec{x}_{0})}{\varepsilon},
\end{equation}
where $\vec{f'}$ doesn't depend on $\vec{h}.$
\end{definition}
A function $\vec{f}$  is said to be hyperholomorphic if $\vec{f}$
is differentiable at every point of $\mathbf{B}$.

\begin{theorem} \label{A}
A function
$\vec{f}(\vec{x})=\sum_{k=0}^{n}\vec{e}_{k}u_{k}(\vec{x})$ is
hyperholomorphic if and only if there exists the function
$\vec{f'}:\mathbf{B}\rightarrow\mathbf{A}$ such that for all
$k=0,1,\ldots,m,$ and $\forall \vec{x}\in \mathbf{B}$
\begin{equation}\label{2}
\vec{e}_{k}\vec{f'}(\vec{x})=\lim_{\varepsilon\rightarrow 0}
\frac{\vec{f}(\vec{x}+\varepsilon\vec{e}_{k})-\vec{f}(\vec{x})}{\varepsilon},
\end{equation}
where $\vec{f'}$  doesn't depend on $\vec{e}_{k}.$
\end{theorem}

\begin{proof} Suppose that (\ref{2}) is fulfilled, then it is easily verified that
\begin{equation}\label{3}
\begin{array}{ccc}\vspace*{3mm}
\vec{f'}=\lim_{\varepsilon \rightarrow 0}
\frac{\vec{f}(\vec{x}+\varepsilon\vec{e}_{0})-\vec{f}(\vec{x})}{\varepsilon}=
\sum_{k=0}^{n}\vec{e}_{k}\frac{\partial u_{k}}{\partial
x_{0}},\\\vspace*{3mm}
\vec{e}_{1}\vec{f'}=\lim_{\varepsilon\rightarrow 0}
\frac{\vec{f}(\vec{x}+\varepsilon\vec{e}_{1})-\vec{f}(\vec{x})}{\varepsilon}=
\sum_{k=0}^{n}\vec{e}_{k}\frac{\partial u_{k}}{\partial x_{1}}=
\vec{e}_{1}\sum_{k=0}^{n}\vec{e}_{k}\frac{\partial u_{k}}{\partial
x_{0}},\\\vspace*{3mm} \vdots\\\vspace*{3mm}
\vec{e}_{m}\vec{f'}=\lim_{\varepsilon\rightarrow 0}
\frac{\vec{f}(\vec{x}+\varepsilon\vec{e}_{m})-\vec{f}(\vec{x})}{\varepsilon}=
\sum_{k=0}^{n}\vec{e}_{k}\frac{\partial u_{k}}{\partial x_{m}}=
\vec{e}_{m}\sum_{k=0}^{n}\vec{e}_{k}\frac{\partial u_{k}}{\partial
x_{0}}.
\end{array}
\end{equation}

Consider $\vec{h}=\sum_{k=0}^{m}h_{k}\vec{e}_{k}.$ It follows from
Eqs.(\ref{3}) that
\[
\begin{array}{ccc}\vspace*{3mm}
h_{0}\vec{f'}= h_{0}\sum_{k=0}^{n}\vec{e}_{k}\frac{\partial
u_{k}}{\partial x_{0}},\\\vspace*{3mm} h_{1}\vec{e}_{1}\vec{f'}=
h_{1}\sum_{k=0}^{n}\vec{e}_{k}\frac{\partial u_{k}}{\partial
x_{1}},\\\vspace*{3mm} \vdots\\\vspace*{3mm}
h_{m}\vec{e}_{m}\vec{f'}=
h_{m}\sum_{k=0}^{n}\vec{e}_{k}\frac{\partial u_{k}}{\partial
x_{m}}.
\end{array}
\]

This implies that
\[
\begin{array}{ccc}\vspace*{3mm}
\vec{h}\vec{f'}= h_{0}\sum_{k=0}^{n}\vec{e}_{k}\frac{\partial
u_{k}}{\partial
x_{0}}+h_{1}\sum_{k=0}^{n}\vec{e}_{k}\frac{\partial
u_{k}}{\partial
x_{1}}+\ldots+h_{m}\sum_{k=0}^{n}\vec{e}_{k}\frac{\partial
u_{k}}{\partial x_{m}}=\\\vspace*{3mm}
\lim\limits_{\varepsilon\rightarrow 0}
\frac{\vec{f}(\vec{x}+\varepsilon\vec{h})-\vec{f}(\vec{x})}{\varepsilon}.
\end{array}
\]

Furthermore, it follows from Eqs.(\ref{3}) that
\[
\begin{array}{ccc}\vspace*{3mm}
h_{0}\sum_{k=0}^{n}\vec{e}_{k}\frac{\partial u_{k}}{\partial
x_{0}}+h_{1}\sum_{k=0}^{n}\vec{e}_{k}\frac{\partial
u_{k}}{\partial
x_{1}}+\ldots+h_{m}\sum_{k=0}^{n}\vec{e}_{k}\frac{\partial
u_{k}}{\partial x_{m}}=\\\vspace*{3mm}
h_{0}\sum_{k=0}^{n}\vec{e}_{k}\frac{\partial u_{k}}{\partial
x_{0}}+h_{1}\vec{e}_{1}\sum_{k=0}^{n}\vec{e}_{k}\frac{\partial
u_{k}}{\partial
x_{0}}+\ldots+h_{m}\vec{e}_{m}\sum_{k=0}^{n}\vec{e}_{k}\frac{\partial
u_{k}}{\partial x_{0}}.
\end{array}
\]

Therefore, for every $\vec{h}\in \mathbf{B}$

$$\vec{h}\sum_{k=0}^{n}\vec{e}_{k}\frac{\partial
u_{k}}{\partial x_{0}}= \lim\limits_{\varepsilon\rightarrow 0}
\frac{\vec{f}(\vec{x}+\varepsilon\vec{h})-\vec{f}(\vec{x})}{\varepsilon}
$$
or
\begin{equation}\label{4}
\vec{f'}=\sum_{k=0}^{n}\vec{e}_{k}\frac{\partial u_{k}}{\partial
x_{0}}
\end{equation}
\end{proof}

By using Eq.(\ref{3}), we have

\begin{equation}\label{5}
\begin{array}{ccc}\vspace*{3mm}
\sum_{k=0}^{n}\vec{e}_{k}\frac{\partial u_{k}}{\partial x_{1}}=
\vec{e}_{1}\sum_{k=0}^{n}\vec{e}_{k}\frac{\partial u_{k}}{\partial
x_{0}},\\\vspace*{3mm} \sum_{k=0}^{n}\vec{e}_{k}\frac{\partial
u_{k}}{\partial x_{2}}=
\vec{e}_{2}\sum_{k=0}^{n}\vec{e}_{k}\frac{\partial u_{k}}{\partial
x_{0}},\\\vspace*{3mm} \vdots\\\vspace*{3mm}
\sum_{k=0}^{n}\vec{e}_{k}\frac{\partial u_{k}}{\partial x_{m}}=
\vec{e}_{m}\sum_{k=0}^{n}\vec{e}_{k}\frac{\partial u_{k}}{\partial
x_{0}}.
\end{array}
\end{equation}

Eqs.(\ref{5}) will be called the Cauchy-Riemann type conditions.
In paper \cite{P} we studied hyperholomorphic functions
$\vec{f}:\mathbf{A}\rightarrow\mathbf{A}$.

\section{Solution of some differential equations}

Let us consider the following partial differential equation

\begin{equation}\label{6}
     \sum_{i_{0}+i_{1}+\ldots +i_{m}=r}C_{i_{0}i_{1}\ldots i_{m}}\frac{\partial^{r}}
     {\partial x^{i_{0}}_{0}\partial x^{i_{1}}_{1}\ldots \partial x^{i_{m}}_{m}}u(x_{0},x_{1},\ldots ,x_{m})=0,
\end{equation}

where $C_{i_{0}i_{1}\ldots i_{m}}\in K$ are constant coefficients.

\begin{theorem}
Suppose that some subspace $\mathbf{B}$ of algebra $\mathbf{A}$
has a basis $\vec{e}_{0}, \vec{e}_{1}, \ldots, \vec{e}_{m}$ such
that
\begin{equation*}
     \sum_{i_{0}+i_{1}+\ldots +i_{m}=r}C_{i_{0}i_{1}\ldots i_{m}}(\vec{e}_{0})^{i_{0}}(\vec{e}_{1})^{i_{1}}\ldots
     (\vec{e}_{m})^{i_{m}}=0,
\end{equation*}
and a function
$\vec{f}(\vec{x})=\sum_{k=0}^{n}\vec{e}_{k}u_{k}(\vec{x})$
$(\vec{f}:\mathbf{B}\rightarrow\mathbf{A})$ is hyperholomorphic.

Then functions $u_{k}(\vec{x})$, $k=0,1,\ldots,n$ are solutions of
Eq.(\ref{6}).
\end{theorem}

\begin{proof}

It follows from the Cauchy-Riemann condition (\ref{5}) that
\begin{equation*}
\frac{\partial \vec{f}}{\partial x_{k}}=\vec{e_{k}}\frac{\partial
\vec{f}}{\partial x_{0}},\ k=0,1,\ldots,m.
\end{equation*}
This implies
\begin{equation*}
\frac{\partial^{i} \vec{f}}{\partial
x_{k}^{i}}=(\vec{e_{k}})^{i}\frac{\partial^{i} \vec{f}}{\partial
x_{0}^{i}}.
\end{equation*}
Therefore, we obtain
\begin{align*}\label{8}
     \sum_{i_{0}+i_{1}+\ldots +i_{m}=r}C_{i_{0}i_{1}\ldots i_{m}}\frac{\partial^{r}}
     {\partial x^{i_{0}}_{0}\partial x^{i_{1}}_{1}\ldots \partial x^{i_{m}}_{m}}\vec{f}(x_{0},x_{1},\ldots ,x_{m})=\\
     \frac{\partial^{r} \vec{f}(x_{0},x_{1},\ldots ,x_{m})}{\partial
x_{0}^{r}}\sum_{i_{0}+i_{1}+\ldots +i_{m}=r}C_{i_{0}i_{1}\ldots
i_{m}}(\vec{e}_{0})^{i_{0}}(\vec{e}_{1})^{i_{1}}\ldots
     (\vec{e}_{m})^{i_{m}}=0,
\end{align*}
\end{proof}

\begin{example}
In case where $\mathbf{A}$ and $\mathbf{B}$ are complex numbers a
differentiable function is holomorphic. Basis $\vec{e}_{0}=1,
\vec{e}_{1}=i$ satisfies the equation
$(\vec{e}_{0})^{2}+(\vec{e}_{1})^{2}=0$. Therefore, the real and
imaginary parts of any holomorphic function satisfy the Laplace
equation
\begin{equation*}
\frac{\partial^{2}}{\partial
x^{2}}u(x,y)+\frac{\partial^{2}}{\partial y^{2}} u(x,y)=0.
\end{equation*}
In paper \cite{K} it is considered algebras which differentiable
functions have harmonic components that satisfy the
three-dimensional Laplace equation.
\end{example}
\bibliographystyle{amsplain}

\end{document}